\documentclass[11pt]{amsart}
\usepackage{amsmath,amsfonts,amssymb,amsthm}
\usepackage{epsfig}
\usepackage{pstricks}
\usepackage{bm}

\newcommand{\Z}{ \mathbb Z}
\newcommand{\N}{ \mathbb N}
\newcommand{\R}{ \mathbb R}

\newcommand{\PP}{ \mathbb P}
\newcommand{\p}[1]{\mathcal{#1}}
\newcommand{\s}[1]{\mathbf{#1}}

\newcommand{\set}[1]{\{#1\}}

\DeclareMathOperator{\supp}{supp}
\DeclareMathOperator{\indeg}{indeg}
\DeclareMathOperator{\outdeg}{outdeg}
\DeclareMathOperator{\lcm}{lcm}
\DeclareMathOperator{\vol}{vol}

\newtheoremstyle{thm}
  {9pt}{9pt}{\itshape}{}{\bfseries}{}{.5em}{}
\theoremstyle{thm}
\newtheorem{thm}{Theorem}[section]
\newtheorem{cor}[thm]{Corollary}
\newtheorem{lemma}[thm]{Lemma}
\newtheorem{prop}[thm]{Proposition}

\newtheorem{alg}[thm]{Algorithm}
\newtheoremstyle{defin}
  {9pt}{9pt}{}{}{\bfseries}{}{.5em}{}
\theoremstyle{defin}

\newtheoremstyle{exm}
  {9pt}{9pt}{}{}{\scshape}{}{.5em}{}
\theoremstyle{exm}
\newtheorem{exm}[thm]{Example}
\newtheorem{rmk}[thm]{Remark}
\newtheoremstyle{proof}
  {}{}{}{}{\itshape}{:}{.5em}{}
\numberwithin{equation}{section}

\title{Geometry and complexity of O'Hara's algorithm}
\author{Matja\v z Konvalinka and Igor Pak}
\date{\today}

\begin{document}

\begin{abstract}
 In this paper we analyze O'Hara's partition bijection. We present three type of results.
 First, we show that O'Hara's bijection can be viewed geometrically as a certain scissor
 congruence type result. Second, we obtain a number of new complexity bounds,
 proving that O'Hara's bijection is efficient in several special cases and mildly
 exponential in general. Finally, we prove that for identities with finite support,
 the map of the O'Hara's bijection can be computed in polynomial time, i.e.\hspace{-0.07cm} much more
 efficiently than by O'Hara's construction.
\end{abstract}

\maketitle

\section{Introduction}

Ever since the pioneer work by Sylvester and his school, there has been a quest to find bijective proofs of many interesting partition identities. Despite remarkable successes in the last century (see~\cite{P3})
and some recent work of both positive and negative nature (see e.g.~\cite{P2,P4}), the problem remains ambiguous and largely unresolved. Much of this stems from the lack of clarity as to what exactly constitutes a bijective proof. Depending on whether one accentuates simplicity, ability to generalize, the time complexity, geometric structure, or asymptotic stability, different answers tend to emerge.

\medskip

In one direction, the subject of partition bijections was revolutionized by Garsia and Milne with their \emph{involution principle}~\cite{GM1,GM2}. This is a combinatorial construction which allows to use a few basic bijections and involutions to build more involved combinatorial maps. As a consequence, one can start with a reasonable analytic proof of a partition identity and trace every step to obtain a (possibly extremely complicated) bijective construction. Garsia and Milne used this route to obtain a long sought bijection proving the Rogers-Ramanujan identities, resolving an old problem in this sense~\cite{GM2}. Unfortunately, this bijection is too complex to be analyzed and has yet to lead to new Rogers-Ramanujan type partition identities.

\medskip

After Garsia-Milne paper, there has been a flurry of activity to obtain synthetic bijections for large classes of partition identities. Most of these bijections did not seem to lead anywhere with one notable exception. Remmel and Gordon found (rather involved) bijective proofs of basic partition identities due to Andrews~\cite{R,G}. The latter are direct extensions of Euler's distinct/odd theorem and have a similar straightforward analytic proof~\cite{A,P3}. Then O'Hara made a surprising discovery that Remmel's and Gordon's bijections can be streamlined to give the same bijective map with a simple construction~\cite{O1,O2}. In fact, O'Hara proved that the resulting bijection is a direct generalization of Glaisher's classical bijection proving Euler's theorem. Moreover, in her  thesis~\cite{O1}, O'Hara showed that her bijection is computationally efficient in certain special cases.
Until now, the reason why O'Hara's bijection has a number of
nice properties distinguishing it from the
other ``involution principle bijections'' remained mysterious.

\medskip

In this paper we obtain results of both positive and negative type.
First, we analyze the complexity of O'Hara's bijection, which we view as a discrete algorithm.
We prove a general result (Theorem~\ref{main1}), which given an \emph{exact} formula for the number of steps of the algorithm in certain cases. From here we show that O'Hara's bijection is computationally
efficient in many special cases. On the other hand, perhaps surprisingly, we prove that
the number of steps can be (mildly) exponential in the worst case (Theorem~\ref{main4} part~\eqref{main4c}). In fact, even when the natural speed-up is applied, the worst complexity does not improve significantly (see Section~\ref{speedy}). This is the first negative result of this kind, proving the analogue of a conjecture that remains open for the Garsia-Milne's ``Rogers-Ramanujan bijection'' (see Subsection~\ref{final1}).

\medskip

Second, we show that O'Hara's bijection has a rich underlying geometry. In a manner similar to that in~\cite{P1,PV}, we view this
bijection as a map between integer points in polytopes which preserves certain linear functionals.
We present an advanced generalization of Andrews's result and of O'Hara's bijection in this geometric setting.
In a special case, the working of the map corresponds to the Euclid algorithm and, more generally,
to terms in the continuing fractions.
Thus one can also think of our generalization as a version of multidimensional continuing fractions.

\medskip

Finally, by combining the geometric and complexity ideas we show that in the finite dimensional case the
map defined by O'Hara's bijection is a solution of an integer linear programming problem.
This implies that the map defined by the bijection can be computed in polynomial time,
i.e.\hspace{-0.07cm} much more efficiently than by O'Hara's bijection. This suggests that perhaps
in the general case O'Hara's bijection can also be sped up to work in polynomial time
(see Subsection~\ref{linprog}).

\medskip

The paper is structured as follows. We start with definitions and notations in Section~\ref{defs}. In Section~\ref{main}, we describe the main results on both geometry and complexity. Proofs of most results are postponed until Section~\ref{proofs}. We conclude with a quick application in Section~\ref{speedy} and final remarks in Section~\ref{final}.

\bigskip

\section{Definitions and background} \label{defs}

\subsection{Andrews's theorem}

A \emph{partition} $\lambda$ is an integer sequence $(\lambda_1,\lambda_2,\ldots,\lambda_\ell)$ such that $\lambda_1 \geq \lambda_2 \geq \ldots \geq \lambda_\ell > 0$, where the integers $\lambda_i$ are called the \emph{parts} of the partition. The sum $n=\sum_{i=1}^\ell \lambda_i$ is called the \emph{size} of $\lambda$, denoted $|\lambda|$; in this case we say that $\lambda$ is a partition of $n$, and write $\lambda \vdash n$. We can also write $\lambda=1^{m_1}2^{m_2}\cdots$, where $m_i=m_i(\lambda)$ is the number of parts of $\lambda$ equal to $i$. The \emph{support} of $\lambda=1^{m_1}2^{m_2}\cdots$ is the set $\set{i \colon m_i > 0}$. The set of all positive integers will be denoted by $\PP$.

\medskip

Denote the set of all partitions by $\p P$ and the set of all partitions of $n$ by $\p P_n$. The number of partitions of $n$ is given by Euler's formula
$$\sum_{\lambda \in \p P} t^{|\lambda|} = \sum_{n=0}^\infty |\p P_n| t^n = \prod_{i=1}^\infty \frac 1{1-t^i}.$$

For a sequence $\overline a = (a_1,a_2,\ldots)$ with $a_i \in \PP \cup \set{\infty}$, define $\p A$ to be the set of partitions $\lambda$ with $m_i(\lambda) < a_i$ for all $i$; write $\p A_n = \p A \cap \p P_n$. Denote by $\supp (\overline a) = \set{i \colon a_i < \infty}$ the \emph{support} of the sequence $a$.

\medskip

Let $\overline a = (a_1,a_2,\ldots)$ and $\overline b = (b_1,b_2,\ldots)$. We say that $\overline a$ and $\overline b$ are \emph{$\varphi$-equivalent}, $\overline a \sim_\varphi \overline b$, if $\varphi$ is a bijection $\supp(\overline a) \to \supp(\overline b)$ such that $i a_i = \varphi(i) b_{\varphi(i)}$ for all $i$. If $\overline a \sim_\varphi\overline b$ for some $\varphi$, we say that $\overline a$ and $\overline b$ are equivalent, and write $\overline a \sim \overline b$.

\begin{thm}[Andrews]
 If $\overline a \sim \overline b$, then $|\p A_n| = |\p B_n|$ for all $n$.
\end{thm}
\begin{proof}
 We use the notation $t^\infty = 0$. Clearly,
 $$\sum_{n = 0}^\infty |\p A_n| t^n = \prod_{i=1}^\infty \frac{1-t^{ia_i}}{1-t^i} = \prod_{j=1}^\infty \frac{1-t^{jb_j}}{1-t^j} = \sum_{n = 0}^\infty |\p B_n| t^n,$$
 which means that $|\p A_n| = |\p B_n|$.
\end{proof}

Consider the classical Euler's theorem on partitions into distinct and odd parts. For $\overline a = (2,2,\ldots)$ and $\overline b = (\infty,1,\infty,1,\ldots)$, $\p A_n$ is the set of all partitions of $n$ into distinct parts, and $\p B_n$ is the set of partitions of $n$ into odd parts. The bijection $i \mapsto 2i$ between $\supp(\overline a)=\PP$ and $\supp(\overline b)=2\PP$ satisfies $ia_i = \varphi(i) b_{\varphi(i)}$, so $\overline a \sim_\varphi \overline b$ and $|\p A_n| = |\p B_n|$. Throughout the paper, we refer to this example as the \emph{distinct/odd case}.

\bigskip \subsection{O'Hara's algorithm}

The analytic proof of Andrews's theorem shown above does not give an explicit bijection $\p A_n \to \p B_n$. Such a bijection is, by Theorem~\ref{defs2}, given by the following algorithm.

%
%

\begin{alg}\label{defs4}
 \emph{(O'Hara's algorithm on partitions)}
 \begin{itemize}
  \item[] \emph{\texttt{Fix:} sequences $\overline a \sim_\varphi \overline b$}
  \item[] \emph{\texttt{Input:} $\lambda \in \p A$}
  \item[] \emph{\texttt{Set:} $\mu \leftarrow \lambda$}
  \item[] \emph{\texttt{While:} $\mu$ contains more than $b_j$ copies of $j$ for some $j$}
  \begin{itemize}
   \item[] \emph{\texttt{Do:} remove $b_j$ copies of $j$ from $\mu$, add $a_i$ copies of $i$ to $\mu$, where $\varphi(i)=j$}
  \end{itemize}
  \item[] \emph{\texttt{Output:} $\psi(\lambda) \leftarrow \mu$}
 \end{itemize}
\end{alg}

\begin{thm}[O'Hara] \label{defs2}
 Algorithm~\ref{defs4} stops after a finite number of steps. The resulting partition $\psi(\lambda) \in \p B$ is independent of the order of the parts removed and defines a size-preserving bijection $\p A \to \p B$.
\end{thm}

Later on (see Subsection~\ref{proofohara}) we deduce O'Hara's theorem from our generalization (Theorem~\ref{main1}).

\medskip

Denote by $L_\varphi(\lambda)$ the number of steps O'Hara's algorithm takes to compute $\psi(\lambda)$, and by $\p L_\varphi(n)$ the maximum value of $L_\varphi(\lambda)$ over all $\lambda \vdash n$.

\begin{exm}
 In the distinct/odd case, O'Hara's algorithm gives the inverse of Glaisher's bijection, which maps $\lambda=1^{m_1}3^{m_3}\cdots \in \p B$ to the partition $\mu \in \p A$ which contains $i 2^j$ if and only if $m_i$ has a $1$ in the $j$-th position when written in binary.
\qed \end{exm}

\begin{exm} \label{defs3}
 Let $\overline a = (1,1,4,5,3,1,1,\ldots)$, $\overline b = (1,1,5,3,4,1,1,\ldots)$ and $\varphi(3)=4$, $\varphi(4)=5$, $\varphi(5)=3$, $\varphi(i)=i$ for $i \neq 3,4,5$; observe that $\overline a \sim_\varphi \overline b$. Then O'Hara's algorithm on $\lambda=3^3 4^4 5^2$ runs as follows:
 $$\begin{array}{cccccccccc}
  & \mathbf{3^3 4^4 5^2} & \to & \scriptstyle 3^7 4^1 5^2 & \to & \scriptstyle 3^2 4^1 5^5 & \to & \scriptstyle 3^2 4^6 5^1 & \to & \scriptstyle 3^6 4^3 5^1 \\
  \to & \scriptstyle 3^{10} 4^0 5^1 & \to & \scriptstyle 3^5 4^0 5^4 & \to & \scriptstyle 3^0 4^0 5^7 & \to & \scriptstyle 3^0 4^5 5^3 & \to &  \mathbf{3^4 4^2 5^3}
 \end{array}$$
 We have $L_\varphi(\lambda) = \p L_\varphi(35) = 9$.
\qed \end{exm}

\begin{exm}
 Take $\overline a = (2,2,1,2,2,1,\ldots)$ and $\overline b=(3,1,3,1,\ldots)$. Here $\p A$ is the set of partitions into distinct parts $\equiv \pm 1$ mod $3$, and $\p B$ is the set of partitions into odd parts, none appearing more than twice. Define $\varphi \colon \PP \to \PP$ as follows:
 \begin{equation} \label{defs1}
  \varphi(i) = \left\{ \begin{array}{cl} i & \mbox{ if } i \mbox{ is divisible by 6} \\ i/3 & \mbox{ if } i  \mbox{ is divisible by 3, but not by 2}\\ 2i & \mbox{ if } i \mbox{ is not divisible by 3} \end{array}\right..
 \end{equation}
 Clearly, $\overline a \sim_\varphi \overline b$. O'Hara's algorithm on $1^1 2^1 8^1 10^1 14^1 20^1$ runs as follows:
 $$\scriptstyle
 \begin{array}{ccccccccccc}
  &  \mathbf{1^1 2^1 8^1 10^1 14^1 20^1} & \to & \scriptstyle 1^1 2^1 8^1 10^3 14^1 & \to & \scriptstyle 1^1 2^1 7^2 8^1 10^3 & \to & \scriptstyle 1^1 2^1 5^2 7^2 8^1 10^2 \\
  \to & \scriptstyle 1^1 2^1 5^4 7^2 8^1 10^1 & \to & \scriptstyle 1^1 2^1 5^6 7^2 8^1 & \to & \scriptstyle 1^1 2^1 4^2 5^6 7^2 & \to & \scriptstyle 1^1 2^3 4^1 5^6 7^2 \\
  \to & \scriptstyle 1^1 2^5 5^6 7^2 & \to & \scriptstyle  1^3 2^4 5^6 7^2 & \to & \scriptstyle  1^5 2^3 5^6 7^2 & \to & \scriptstyle 1^7 2^2 5^6 7^2 \\
  \to & \scriptstyle 1^9 2^1 5^6 7^2 & \to & \scriptstyle 1^{11} 5^6 7^2 & \to & \scriptstyle 1^{11} 5^3 7^2 15^1 & \to & \scriptstyle 1^{11} 7^2 15^2 \\
  \to & \scriptstyle 1^8 3^1 7^2 15^2 & \to & \scriptstyle 1^5 3^2 7^2 15^2 & \to & \scriptstyle 1^2 3^3 7^2 15^2 & \to &  \mathbf{1^2 7^2 9^1 15^2}
 \end{array}$$
 The bijection $\psi$ is similar in spirit to Glaisher's bijection: given $\lambda=1^{m_1}2^{m_2}4^{m_4}5^{m_5}\cdots \in \p A$ and $j \in \PP$, the number of copies of part $2j-1$ in $\psi(\lambda)$ is equal to the $k$-th digit in the ternary expansion of $l$, where $k$ is the highest power of $3$ dividing $2j-1$, $2j-1=3^k r$, and $l = \sum_i 2^im_{r2^i}$.
\qed \end{exm}

\bigskip \subsection{Equivalent sequences and graphs} \label{graphs}

Choose equivalent sequences $\overline a$, $\overline b$. Define a directed graph $G_\varphi$ on $\supp (\overline a) \cup \supp (\overline b)$ by drawing an edge from $i$ to $j$ if $\varphi(j) = i$; an arrow from $i$ to $j$ therefore means that O'Hara's algorithm simulta\-neously removes copies of $i$ and adds copies of $j$. Each vertex $v$ has $\indeg v \leq 1$, $\outdeg v \leq 1$ and $\indeg v + \outdeg v \geq 1$. The graph splits into connected components of the following five types:

\begin{enumerate}
 \renewcommand{\theenumi}{\roman{enumi}}
 \item \label{defs5a} cycles of length $m \geq 1$;
 \item \label{defs5b} paths of length $m \geq 2$;
 \item \label{defs5c} infinite paths with a starting point, but without an ending point;
 \item \label{defs5d} infinite paths with an ending point, but without a starting point;
 \item \label{defs5e} infinite paths without a starting point or an ending point.
\end{enumerate}

\begin{figure}[ht]
  \begin{center}
   \input{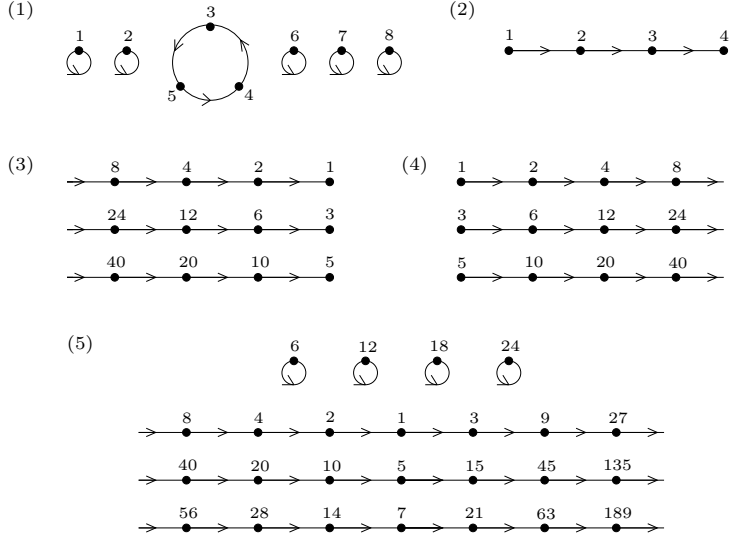}
  \end{center}
  \caption{Examples of graphs $G_\varphi$.}
  \label{fig1}
 \end{figure}

\begin{exm}
 Figure~\ref{fig1} shows portions of graphs $G_\varphi$ for certain $\varphi$:
 \begin{enumerate}
  \item $\overline a = (1,1,4,5,3,1,1,\ldots)$, $\overline b = (1,1,5,3,4,1,1,\ldots)$, $\varphi(3)=4$, $\varphi(4)=5$, $\varphi(5)=3$, $\varphi(i)=i$ for $i \neq 3,4,5$; components of $G_\varphi$ are of type~\eqref{defs5a};
  \item $\overline a = (\infty,1,2,3,\infty,\infty,\infty,\ldots)$, $\overline b=(2,3,4,\infty,\infty,\infty,\infty,\ldots)$, $\varphi(2)=1$, $\varphi(3)=2$, $\varphi(4)=3$; $G_\varphi$ is of type~\eqref{defs5b};
  \item the distinct/odd case: $\overline a = (2,2,\ldots)$, $\overline b = (\infty,1,\infty,1,\ldots)$, $\varphi(i) = 2i$; components of $G_\varphi$ are of type~\eqref{defs5c};
  \item the odd/distinct case: $\overline a = (\infty,1,\infty,1,\ldots)$, $\overline b = (2,2,\ldots)$, $\varphi(i) = i/2$; components of $G_\varphi$ are of type~\eqref{defs5d};
  \item $\overline a = (2,2,1,2,2,1,\ldots)$ and $\overline b=(3,1,3,1,\ldots)$, $\varphi$ given by~\eqref{defs1}; components of $G_\varphi$ are of types~\eqref{defs5a} and~\eqref{defs5e}. \qed
 \end{enumerate}
\end{exm}

\bigskip \subsection{Scissor-congruence and $\Pi$-congruence}

We say that convex polytopes $A,B$ in $\R^m$ are \emph{congruent}, write $A \simeq B$, if $B$ can be obtained from $A$ by rotation and translation. For convex polytopes $P,Q \subset \R^m$, we say that they are \emph{scissor-congruent} if $P$ can be cut into finitely many polytopes which can be rearranged and assembled into $Q$, i.e.\hspace{-0.07cm} if $P$ and $Q$ are the disjoint union of congruent polytopes: $P=\cup_{i=1}^n P_i$, $Q=\cup_{i=1}^n Q_i$, $P_i \simeq Q_i$.

\medskip

Let $\pi$ be a linear functional on $\R^m$. If $Q_i$ can be obtained from $P_i$ by a translation by a vector in the hyperplane $\p H= \set{\s x \in \R^m \colon \pi(\s x) = 0}$, we say that $P$ and $Q$ are \emph{$\pi$-congruent}. If $P$ and $Q$ are $\pi$-congruent for some linear functional $\pi$, we say that they are \emph{$\Pi$-congruent}.

\medskip

If $P$ can be cut into countably many polytopes which can be translated by a vector in the hyperplane $\p H = \set{\s x \in \R^m \colon \pi(\s x) = 0}$ and assembled into $Q$, we say that $P$ and $Q$ are \emph{approximately $\pi$-congruent}. We say that they are \emph{approximately $\Pi$-congruent} if they are approximately $\pi$-congruent for some linear functional $\pi$. If $P$ and $Q$ are approximately $\pi$-congruent, there exist, for every $\varepsilon>0$, $\pi$-congruent polytopes $P_\varepsilon \subseteq P$ and $Q_\varepsilon \subseteq Q$, such that $\vol(P \setminus P_\varepsilon) < \varepsilon$ and $\vol(Q \setminus Q_\varepsilon) < \varepsilon$.

\medskip

Finally, let $\s R(a_1,\ldots,a_m)=[0,a_1) \times \cdots \times [0,a_m)$ be a box in $\R^m$, and let $R(a_1,\ldots,a_m) = \s R(a_1,\ldots,a_m) \cap \Z^m$ be the set of its integer points.

\begin{exm}\label{ex-euclid}
 Let $d=2$ and $\pi(x,y)=x+y$. Euclid's algorithm on $(a,b)$ yields a $\pi$-congruence between $\s R(a,b)$ and $\s R(b,a)$: if $b = r_1 a + s_1$ with $0 \leq s_1 < a$, divide $[0,a) \times [0,r_1 a)$ into $r_1$ squares with side $a$, and translate the square $[0,a) \times [i a, (i+1) a)$ by the vector $(i a, - i a)$ to $[i a, (i+1) a) \times [0,a)$. Then write $a = r_2 s_1 + s_2$ with $0 \leq s_2 < s_1$, divide $[0,a) \times [r_1 a, b)$ into $r_2$ squares with side $s_1$, and translate the square $[i s_1, (i+1) s_1) \times [r_1 a, b)$ by the vector $(r_1 a - i s_1,i s_1 - r_1 a)$ to $[r_1 a, b) \times [i s_1, (i+1) s_1)$. Continue until the remainder $s_i$ is equal to $0$. The first drawing of Figure~\ref{fig5} gives an example.
 
 \smallskip 
 
 The second drawing shows that boxes $\s R(12,8)$ and $\s R(32,3)$ are $\pi$-congruent for $\pi(x,y)=x+4y$. Finally, in Figure~\ref{fig2} we give a $\pi$-congruence between $\s R(4,5,3)$ and $\s R(5,3,4)$ for $\pi(x,y,z)=3x+4y+5z$. \qed
 \begin{figure}[ht]
  \begin{center}
   \input{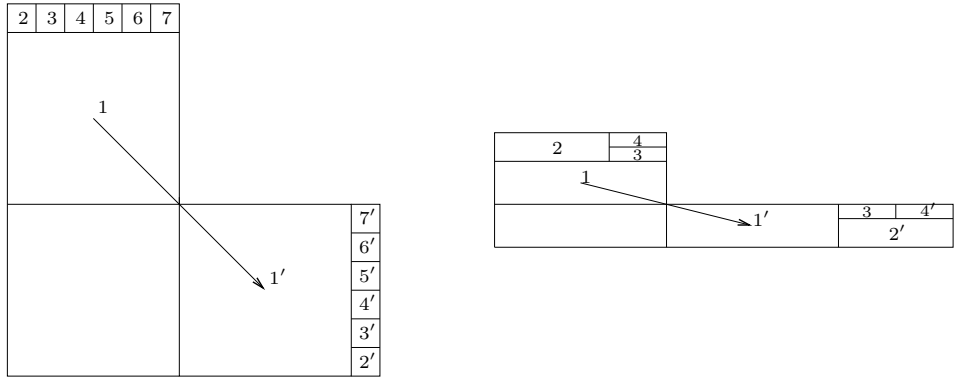}
  \end{center}
  \caption{Two $\Pi$-congruences.}
  \label{fig5}
 \end{figure}

 \begin{figure}[ht]
  \begin{center}
   \epsfig{file=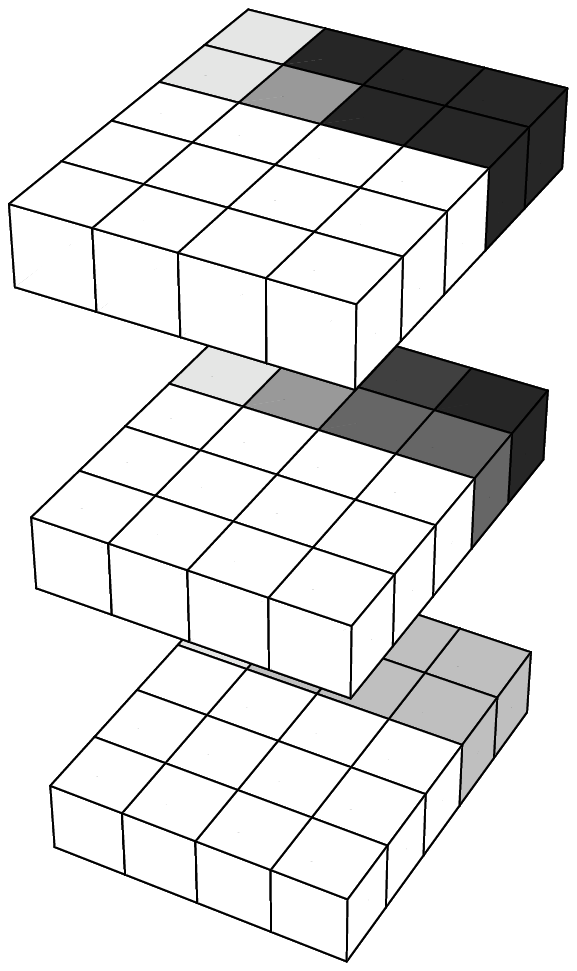,height=3.3cm}
   \begin{minipage}{2cm}
    \begin{center}
     \vspace{-2cm}
     $\stackrel{\bm \psi}{\longrightarrow} \!\!$
     \vspace{2cm}
    \end{center}
   \end{minipage}
   \epsfig{file=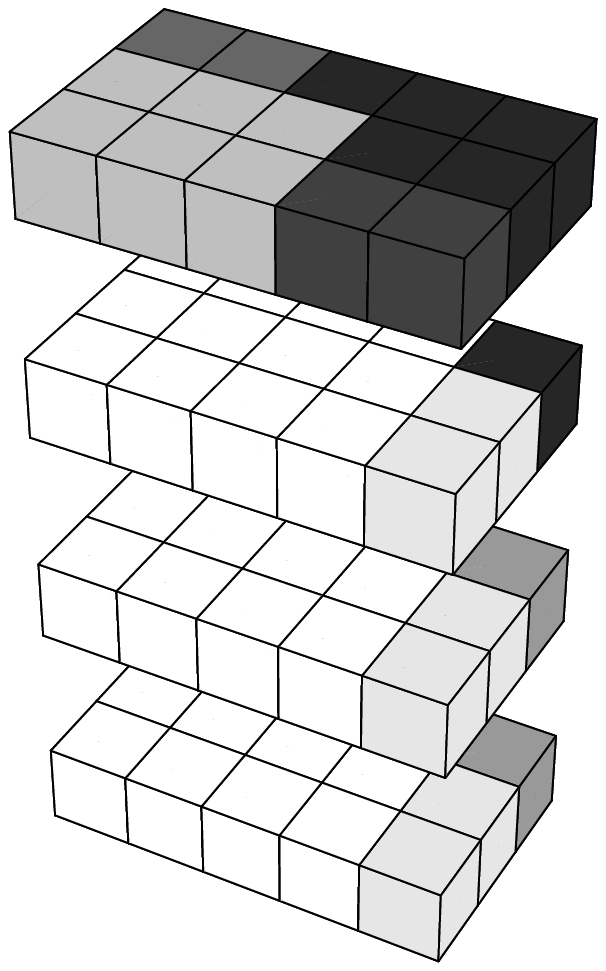,height=3.3cm}
  \end{center}
  \caption{$\pi$-congruence between $\s R(4,5,3)$ and $\s R(5,3,4)$.}
  \label{fig2}
 \end{figure}
\end{exm}

\bigskip

\section{Main results} \label{main}

\subsection{Continuous O'Hara's algorithm and $\Pi$-congruences}

Take the case when $G_\varphi$ is a cycle $i_1 \to i_m \to i_{m-1} \to \ldots \to i_1$. In this case, $\varphi(i_1)=i_2$, $\varphi(i_2)=i_3$, etc. Throughout this section, identify a partition $i_1^{t_1}\cdots i_m^{t_m}$ with the vector $\s t = (t_1,\ldots,t_m)$. By Theorem~\ref{defs2}, O'Hara's algorithm defines a bijection $\psi \colon R(a_1,\ldots,a_m) \to R(b_1,\ldots,b_m)$, where $i_j a_j = i_{j+1} b_{j+1}$ for all $j$. The following algorithm (see also Theorem~\ref{main1}) generalizes $\psi$ to the continuous setting. It gives a bijection ${\bm \psi} \colon \s R(a_1,\ldots,a_m) \to \s R(b_1,\ldots,b_m)$, which is defined also for non-integer $a_j,b_j$. When $a_j,b_j$ are integers, it is an extension of $\psi \colon R(a_1,\ldots,a_m) \to R(b_1,\ldots,b_m)$. As an immediate corollary, we prove that two boxes with rational coordinates and with equal volume are $\Pi$-congruent. We show in Subsection~\ref{proofohara} how we can use Theorem~\ref{main1} to give an alternative proof of Theorem~\ref{defs2}.

\begin{alg}\label{main7}
 \emph{(continuous O'Hara's algorithm)}
 \begin{itemize}
  \item[] \emph{\texttt{Fix:} $\s i = (i_1,\ldots,i_m) \in \R_+^m$}
  \item[] \emph{\phantom{\texttt{Fix:}} $\s a = (a_1,\ldots,a_m) \in \R_+^m$, $\s b = (b_1,\ldots,b_m) \in \R_+^m$ with $i_ja_j = i_{j+1} b_{j+1}$}
  \item[] \emph{\texttt{Input:} $\s t \in \s R(a_1,\ldots,a_m)$}
  \item[] \emph{\texttt{Set:} $\s s \leftarrow \s t$}
  \item[] \emph{\texttt{While:} $\s s$ contains a coordinate $s_j \geq b_j$}
  \begin{itemize}
   \item[] \emph{\texttt{Do:} $s_j \leftarrow s_j-b_j$, $s_{j-1} \leftarrow s_{j-1}+a_{j-1}$}
  \end{itemize}
  \item[] \emph{\texttt{Output:} $\bm \psi(\s t) \leftarrow \s s$}
 \end{itemize}
\end{alg}

It is clear that the algorithm starts with an element of $P = \s R(a_1,\ldots,a_m)$ and, if the while loop terminates, outputs an element of $Q = \s R(b_1,\ldots,b_m)$. It is not obvious, however, that the loop terminates in every case, or that the output $\bm \psi(\s t)$ and the number of steps $\s L_\varphi(\s t)$ depend only on $\s t$, not on the choices made in the while loop.

\begin{thm} \label{main1}
 Algorithm~\ref{main7} has the following properties.
 \begin{enumerate}
  \item \label{main1a} The algorithm stops after a finite number of steps, and the resulting vector $\bm \psi(\s t)$ and the number of steps $\s L_\varphi(\s t)$ are independent of the choices made during the execution of the algorithm.
  \item \label{main1b} The algorithm defines a bijection $\bm \psi \colon P \to Q$ which satisfies $\bm \psi(\s t) - \s t \in \p H$, where $\p H$ is the hyperplane defined by $i_1x_1+\ldots+i_mx_m=0$.
  \item \label{main1d} We have
  $$\s L_\varphi(\s t + \s t') \geq \s L_\varphi(\s t) + \s L_\varphi(\s t') \mbox{ for every } \s t,\s t',\s t+\s t' \in P.$$
  In particular, $\s L_\varphi(\s t') \leq \s L_\varphi(\s t)$ if $\s t' \leq \s t$.
  \item \label{main1c} Let $\s t, \s t' \in P$, $\s s = \bm \psi(\s t)$, with $t_j \leq t_j' < t_j + \varepsilon_j$, where $\varepsilon_j = b_j - s_j$. Then
  $$\bm \psi(\s t')-\s t'=\bm \psi(\s t)-\s t \quad \mbox{and} \quad \s L_\varphi(\s t')=\s L_\varphi(\s t).$$
  \item \label{main1e} For all $\s a, \s b \in \Z_+^m$, we have
  $$\max_{\s t \in P} \s L_\varphi(\s t) \, = \, \lcm(c_1,\ldots,c_m) \cdot \left( \frac 1 {c_1} + \ldots + \frac 1 {c_m} \right) - m,$$
  where $c_j = a_{1}\cdots a_{{j-1}}b_{j}\cdots b_{{m-1}}$.
 \end{enumerate}
\end{thm}

The proof of the theorem is given in Subsection~\ref{proofmain1}.

\medskip

We call boxes $P = \s R(a_1,\ldots,a_m)$, $Q=\s R(b_1,\ldots,b_m)$ \emph{relatively rational} if there exists $\lambda$, $\lambda \neq 0$, such that $\lambda a_j \in \Z, \lambda b_j \in \Z$. Clearly, two boxes $P$ and $Q$ with rational side-lengths are relatively rational.

\begin{cor}
 Boxes $P=\s R(a_1,\ldots,a_m)$, $Q=\s R(b_1,\ldots,b_m)$ with equal volume are approximately $\Pi$-congruent. Moreover, when $P$ and $Q$ are relatively rational and have equal volume, they are $\Pi$-congruent.
\end{cor}
\begin{proof}
 For $j=1,\ldots,m$, take $i_j = a_1\cdots a_{j-1}b_{j+1}\cdots b_m$. Clearly $i_j a_j = i_{j+1} b_{j+1}$ for $j=1,\ldots,m-1$, and $a_1\cdots a_m = b_1\cdots b_m$ implies $i_m a_m = i_1 b_1$. Therefore, the numbers $i_j,a_j,b_j$ satisfy the conditions of Algorithm~\ref{main7}. By Theorem~\ref{main1} part~\eqref{main1b}, the algorithm defines a bijection $\bm \psi \colon P \to Q$. Parts~\eqref{main1c} and~\eqref{main1b} of Theorem~\ref{main1} imply that we can cut $P$ into (countably many) smaller boxes, each of which is translated by a vector in the plane $i_1x_1+\ldots+i_mx_m=0$.
 
 \smallskip
 
 If $P$ and $Q$ are relatively rational, we can assume without loss of generality that all $a_j,b_j$ are integers. For any integer vector $\s t$, we have $\bm \psi(\s t')-\s t' = \bm \psi(\s t)-\s t$ and $\s L_\varphi(\s t')=\s L_\varphi(\s t)$ whenever $t_j \leq t_j' < t_j + 1$, so $P$ and $Q$ are divided into a finite number (at most $a_1\cdots a_m$) of boxes.
\end{proof}

\begin{exm}
 Even in the $3$-dimensional case the $\Pi$-congruence defined by the algorithm can be quite complex,
 as the next figure suggests. Here the same shading is used for parallel translations by the same vector.
\qed
 \begin{figure}[ht]
  \begin{center}
   \epsfig{file=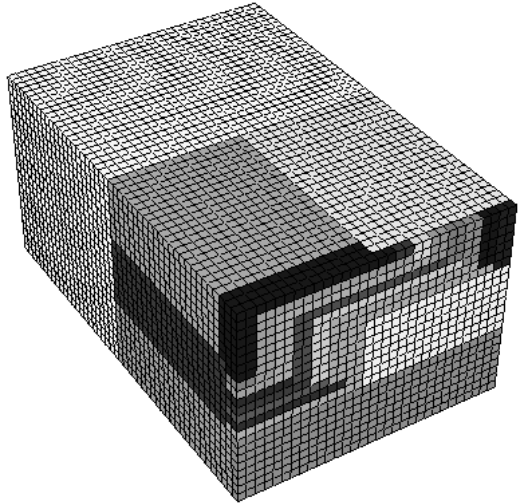,height=7cm}
  \end{center}
  \caption{The decomposition of the box $\s R(31,47,23)$ given by O'Hara's algorithm (only the top, right, and back sides are shown) .}
  \label{fig3}
 \end{figure} \end{exm}

\bigskip
\subsection{Complexity of O'Hara's algorithm}

The complexity of O'Hara's algorithm has been an open problem, with the exception of the elementary distinct/odd case (see~\cite{O1} and Example~\ref{proofscxity1}).

\medskip

It turns out that the complexity depends heavily on the type of the graph $G_\varphi$ defined in Subsection~\ref{graphs}. Part~\eqref{main1e} of Theorem~\ref{main1} gives the maximum number of steps that O'Hara's algorithm takes when $G_\varphi$ is a cycle. The following lemma gives an estimate for $\p L_\varphi(n)$ when $G_\varphi$ is a path.

\begin{lemma} \label{main2}
 Let $G_\varphi$ be a finite or infinite path on $\p I \subseteq \PP$. Then $\p L_\varphi(n) \leq n (\log n + 1)$. Moreover, if
 $$D = \sum_{i \in \p I} \frac 1 {i a_i} = \sum_{j \in \p I} \frac 1 {j b_j} < \infty,$$
 then $\p L_\varphi(n) \leq D n$. 
\end{lemma}

This lemma and the other results in this subsection are proved in Section~\ref{proofs}. Here and throughout the paper, by $\log n$ we mean the natural logarithm of $n$.

\medskip

We combine these estimates to prove Theorem~\ref{main3}, the second main result of this paper, which gives upper bounds for the complexity of the algorithm in the general case.

\begin{thm} \label{main3}
 Let $\overline a,\overline b$ be $\varphi$-equivalent sequences.
 \begin{enumerate}
  \item \label{main3a} If $G_\varphi$ has only a finite number of cycles of length $> 2$, then $\p L_\varphi(n) = O(n \log n)$, and the constants implied by the $O$-notation are universal.
  \item \label{main3b} If $G_\varphi$ has only a finite number of cycles of length $> m$ for some $m > 2$, then $\p L_\varphi(n) = O(n^{m-1})$, and the constants implied by the $O$-notation depend only on $m$.
 \end{enumerate}
\end{thm}

The following theorem gives the corresponding lower bound on the worst case complexity. It shows that the estimates of Theorem~\ref{main3} are close to being sharp.

\begin{thm} \label{main4}
 There exist $\varphi$-equivalent sequences $\overline a$ and $\overline b$, such that:
 \begin{enumerate}
  \item \label{main4a} $G_\varphi$ is a path and $\p L_\varphi(n) = \Omega(n \log \log n)$;
  \item \label{main4b} $G_\varphi$ contains only cycles of length $\leq m$ and $\p L_\varphi(n) = \Omega(n^{m-1-\varepsilon})$ for every $\varepsilon>0$;
  \item \label{main4c} $\p L_\varphi(n) = \exp \Omega(\sqrt[3]{n})$.
 \end{enumerate}
\end{thm}

In other words, depending on the type of the graph, we have nearly matching upper and
lower bounds on $\p L_\varphi(n)$. For example, for an $m$-cycle, Theorem~\ref{main3}
shows that $\p L_\varphi(n)$ is $O(n^{m-1})$, while Theorem~\ref{main4} shows that it
is $\Omega(n^{m-1-\varepsilon})$ for every $\varepsilon>0$. Similarly, part~\eqref{main4c}
shows that O'Hara's algorithm can be very slow in general since the total number of
partitions of~$n$ is asymptotically  $\exp \Theta(\sqrt n)$.

\begin{exm} \label{proofscxity1}
 In the distinct/odd case, the graph $G_\varphi$ is composed of infinite paths
 $$\ldots \to 8 j \to 4 j \to 2j \to j \quad \mbox{for each odd } j.$$
 A partition $\lambda \vdash n$ can be broken up into partitions $\lambda^{(j)} \vdash n_j$ such that the support of $\lambda^{(j)}$ is contained in $\set{(2j-1)2^k \colon k \in \N}$. We have
 $$\sum_{k=0}^\infty \frac 1 {2^{k+1} (2j-1)} = \frac 1 {2j-1},$$
 and Lemma~\ref{main2} implies that O'Hara's algorithm takes at most $n_j/(2j-1)$ steps to compute $\psi(\lambda^{(j)})$. This implies that
 $$L_\varphi(\lambda) \, \leq \, \sum_{j=1}^\infty \frac{n_j}{2j-1} \, \leq \, \sum_{j=1}^n n_j \, = \, n.$$
 In other words, O'Hara's algorithm takes at most $n$ steps to compute $\psi(\lambda)$. This bound is (almost) sharp since the algorithm takes $2^k-1$ steps to compute $\psi(2^k)=1^{2^k}$.
\qed \end{exm}

\bigskip \subsection{O'Hara's algorithm as an integer linear programming problem} \label{linprog}

Let us now give a new description of O'Hara's algorithm.

\begin{prop} \label{main6}
 Let $\s i, \s a, \s b \in$ be as above such that $i_j a_j = i_{j+1} b_{j+1}$ for $j=1,\ldots,m$. Fix a vector $\s t \in \s R(a_1,\ldots,a_m)$. Then $\s s = \bm \psi(\s t)$ satisfies the following:
 $$\s s = \s t + A \s k,$$
 where
 $$A = \begin{pmatrix} -b_1 & a_1 & 0 & \cdots & 0 \\ 0 & -b_2 & a_2 & \cdots & 0 \\ 0 & 0 & -b_3 & \cdots & 0 \\ \vdots & \vdots & \vdots & \ddots & \vdots \\ a_m & 0 & 0 & \cdots & -b_m\end{pmatrix}$$
 and $\s k=(k_1,\ldots,k_m)$ is the unique vector minimizing
 $$k_1+\ldots+k_m$$
 with constraints
 $$\s k \in \Z^m, \qquad \s k \geq \s 0, \qquad A \s k \geq - \s t, \qquad A \s k \leq \s b - \s 1 - \s t.$$
\end{prop}

The proposition is proved in Subsection~\ref{proofmain6}. The advantage of this approach is that one can use standard integer linear programming results to speed up the computation of $\bm \psi(\s t)$.

\begin{thm} \label{linprog1}
 For every $m \geq 1$, there exists a deterministic algorithm which computes the continuous O'Hara's bijection $\bm \psi$ in polynomial time, for all integer vectors $\s i,\s a,\s b$ as above.
\end{thm}
\begin{proof}
 It is well known that for a bounded dimension $m$, there exists an algorithm for solving integer linear programming problem $A\s x \leq \s b$, for which the number of steps is bounded by a polynomial in the logarithm of the largest entry of $A$, $\s b$ for integer $A$ and $\s b$ (see e.g.\hspace{-0.07cm}~\cite[Corollary 18.7b]{S}). By Proposition~\ref{main6}, this implies the result.
\end{proof}

Theorem~\ref{linprog1} can be used to obtain a significant speed-up of (the usual) O'Hara's algorithm, in the case when $G_\varphi$ contains only cycles of bounded length. Namely, we obtain the following result.

\begin{thm} \label{linprog2}
 Let $\overline a \sim_\varphi \overline b$. If the lengths of cycles of $G_\varphi$ are bounded, there exists a deterministic algorithm which computes $\psi(\lambda)$ in $O(n \log n)$ steps for $\lambda \in \p A_n$.
\end{thm}
\begin{proof}
 Without loss of generality, the support of $\lambda \in \p A_n$ is contained in one of the connected components of $G_{\varphi}$. If this connected component is a path, O'Hara's algorithm takes $O(n \log n)$ steps by Lemma~\ref{main2}. If it is a cycle of length $m$, the algorithm described in the previous theorem takes $O(\log^{c}n)$ steps for some $c$, and obviously the $O(n \log n)$ term dominates.
\end{proof}

\begin{rmk}
 Let us note that the inner workings of the algorithms in Theorem~\ref{linprog1} and Theorem~\ref{linprog2} have a geometric rather than combinatorial nature, and are very different from those of O'Hara's algorithm. However, both kinds of algorithms, when applied to the same input, have the same output, which means that they produce the same partition bijection.
\end{rmk}

\bigskip

\section{Proofs and examples} \label{proofs}

\subsection{Proof of Theorem~\ref{main1}} \label{proofmain1}

Throughout the section, indices are taken modulo $m$.

\begin{lemma}
 Take a vector $\s t \in P$, choose $\varepsilon_j < a_j-t_j$, and do the algorithm on $\s t$: denote the vectors we get by $\s s^0 = \s t$, $\s s^1$, $\s s^2$, etc. Then each box of size $\varepsilon_1 \times \cdots \times \varepsilon_m$ contains at most one of $\s s^i$.
\end{lemma}
\begin{proof}
  Assume that we have $|s^{k'}_j - s^k_j| \leq \varepsilon_j$ for all $j$ for some $0 \leq k < k'$, i.e.\hspace{-0.07cm} that we hit an $\varepsilon_1 \times \cdots \times \varepsilon_m$ box twice. Say that in the course of getting from $\s t$ to $\s s^k$ (respectively $\s s^{k'}$, respectively $\s s^{k'-k}$), we subtracted $b_j$ from the $j$-th coordinate and added $a_j$ to the $(j-1)$-th coordinate $k_j$ times (respectively $k_j'$ times, respectively $k_j''$ times). Clearly we have $\sum k_j = k$, $\sum k_j'=k'$ and $\sum k_j'' = k'-k$. Furthermore, the equations $s^k_j = t_j - b_j k_j + a_j k_{j+1}$, $s^{k'}_j = t_j - b_j k_j' + a_j k_{j+1}'$ and $s^{k'-k}_j = t_j - b_j k_j'' + a_j k_{j+1}''$, which hold for all $j$, can be written as
 \begin{eqnarray*}
  \s s^k & = & \s t + A \s k, \\
  \s s^{k'} & = & \s t + A \s k', \\
  \s s^{k'-k} & = & \s t + A \s k'',
 \end{eqnarray*}
 where
 $$A = \begin{pmatrix} -b_{1} & a_{1} & 0 & \cdots & 0 \\ 0 & -b_{2} & a_{2} & \cdots & 0 \\ 0 & 0 & -b_{3} & \cdots & 0 \\ \vdots & \vdots & \vdots & \ddots & \vdots \\ a_{m} & 0 & 0 & \cdots & -b_{m}\end{pmatrix}.$$
 There are two cases to consider: $\s k'' = \s k' - \s k$ and $\s k'' \neq \s k' - \s k$. We obtain a contradiction in each case. Assume first that $\s k'' = \s k' - \s k$. Then we have
 $$\s s^{k'-k} - \s t = \s t + A(\s k' - \s k) - \s t = \s s^{k'} - \s s^k.$$
 In particular, we have $|s^{k'-k}_j - t_j| \leq \varepsilon_j$ and $s^{k'-k}_j < a_j$ for all $j$. On the other hand, $\s s^{k'-k}$ was obtained from $\s s^{k'-k-1}$ by choosing $j$ with $s^{k'-k-1}_j \geq b_j$ and then taking $s^{k'-k}_j = s^{k'-k-1}_j - b_j$ and $s^{k'-k}_{j-1} = s^{k'-k-1}_{j-1}+a_{j-1}$. Therefore $s^{k'-k}_{j-1} \geq a_{j-1}$, a contradiction.
 
 \smallskip
 
 Suppose now that $\s k'' \neq \s k' - \s k$. Since $\sum k_j'' = \sum (k_j'-k_j)$, there is a $j$ so that $k_j''\geq k_j'-k_j$, $k_{j+1}'' < k_{j+1}'-k_{j+1}$. We have
 $$\s s^{k'-k} = \s t + A \s k'' = \s t + A(\s k''- \s k'+ \s k) + \s s^{k'} - \s s^k.$$
 This implies
 $$0 \leq s_j^{k'-k} \leq t_j - b_j \underbrace{(k_j''-k_j'+k_j)}_{\geq 0} + a_j\underbrace{(k_{j+1}''-k_{j+1}'+k_{j+1})}_{\leq -1} + \varepsilon_j \leq t_j - a_j + \varepsilon_j < 0,$$
 a contradiction. Therefore we can never hit an $\varepsilon_1 \times \cdots \times \varepsilon_m$ box twice, which completes the proof.
\end{proof}

Now we are ready to prove all parts of Theorem~\ref{main1}.

\medskip

\eqref{main1a} Clearly we have $i_1 s^k_1 + \ldots + i_m s^k_m = i_1 t_1 + \ldots + i_m t_m$ and $s_j^k \geq 0$ for all $j$ and $k$. Since the set $\set{\s x \in \R^m \colon i_1x_1+\ldots+i_mx_m = C, x_j\geq 0}$ can be covered by a finite number of boxes of size $\varepsilon_1 \times \cdots \times \varepsilon_m$, the number of steps is finite.

\medskip

Let us prove that the output vector $\s s \in Q$ does not depend on the choices we make in the course of O'Hara's algorithm. Assume that when we run the algorithm twice, we obtain vectors $\s s$ and $\s s'$. Denote by $k=|\s k|=k_1+\ldots+k_m$ the number of steps to obtain $\s s$, and by $k'=|\s k'|$ the number of steps to obtain $\s s'$, and assume that $k \leq k'$. Denote by $\s s''$ the vector we obtain after $k$ steps in the second run of the algorithm. Think of $\s s''$ as being on the path from $\s t$ to $\s s'$. Write $\s s = \s t + A \s k$, $\s s'' = \s t + A \s k''$, where $|\s k| = |\s k''| = k$. Then $\s s'' = \s s + A (\s k'' - \s k)$, and $k_j'' > k_j$, $k_{j+1}'' \leq k_{j+1}$ would imply
$$0 \leq s''_j = s_j - b_j \underbrace{(k_j''-k_j)}_{\geq 1} + a_j \underbrace{(k_{j+1}''-k_{j+1})}_{\leq 0} < b_j-b_j = 0,$$
which is a contradiction. Therefore $\s k'' = \s k$ and $\s s''=\s s$, therefore also $k'=k$ and $\s s'=\s s$.

\medskip

\eqref{main1b} Let us construct an explicit inverse map $\bm{\psi}^{-1}$. Denote by $\tau$ the flip $\tau \colon (x_1,\ldots,x_m) \mapsto (x_m,\ldots,x_1)$. Consider a map $\bm {\psi'} \colon \s R(b_m,\ldots,b_1) \to \s R(a_m,\ldots,a_1)$ obtained by replacing $(i_1,\ldots,i_m)$, $(a_1,\ldots,a_m)$ and $(b_1,\ldots,b_m)$ with $(i_m,\ldots,i_1)$, $(b_m,\ldots,b_1)$ and $(a_m,\ldots,a_1)$, respectively. Observe that
$$\bm {\psi}^{-1} := \tau \circ \bm {\psi'} \circ \tau \, \colon Q \, \to P$$
is the inverse map of $\bm \psi$. Therefore, $\bm \psi$ is one-to-one.

\medskip

It remains to check that $\bm{\psi}$ satisfies $\bm \psi(\s t) - \s t \in \p H$, where $\p H$ is the hyperplane defined by $i_1x_1+\ldots+i_mx_m=0$. Note that the columns of $A$ lie in~$\p H$. Because $\bm \psi(\s t) - \s t = A \s k$ for some $\s k$, $\bm \psi(\s t) - \s t$ lies in~$\p H$.

\medskip

\eqref{main1d} Let $\s s^0 = \s t,\s s^1,\s s^2,\ldots,\s s^k = \bm{\psi}(\s t)$, $k = \s L_\varphi(\s t)$, be the intermediate steps of the algorithm which computes $\bm{\psi}(\s t)$. Similarly, let $\s r^0 = \s t',\s r^1,\s r^2,\ldots,\s r^{k'} = \bm{\psi}(\s t')$, $k' = \s L_\varphi(\s t')$, be the intermediate steps of the algorithm which computes $\bm{\psi}(\s t)$. Every vector in the sequence $\s t + \s t',\s s^1 + \s t',\ldots,\s s^k + \s t',\s s^k + \s r^1,\ldots,\s s^k + \s r^{k'-1}$ has at least one of the coordinates $\geq b_j$, so O'Hara's algorithm takes at least $k+k'=\s L_\varphi(\s t)+\s L_\varphi(\s t')$ steps to compute $\bm \psi(\s t+\s t')$.

\medskip

\eqref{main1c} Let $\s t \in P$, $\s s = \bm \psi(\s t)$ and $\varepsilon_j = b_j - s_j$. For every $\s t'$ with $t_j \leq t_j' < t_j + \varepsilon_j$, make the steps of the algorithm which inputs $\s t'$ the same as the one which inputs~$\s t$. Let $\s s^i$ and $\s r^i$ be as in part~\eqref{main1d}. Then the vector $\s r^i$ satisfies $\s r^i = \s s^i + (\s t' - \s t)$. For $k=\s L_\varphi(\s t)$, we have $\s s^k=\s s$ and $\s r^k=\s s + (\s t' - \s t) \in Q$. Therefore, $\s r^k = \bm \psi(\s t')$, $\psi(\s t') = \bm \psi(\s t) + \s t' - \s t$ and $\s L_\varphi(\s t')=\s L_\varphi(\s t)$, as desired.

\medskip

\eqref{main1e} Assume that $a_j$ and $b_j$ are integers for all $j=1,\ldots,m$. Let $\s t \in P$ and assume that O'Hara's algorithm with input $\s t$ subtracts $b_j$ from the $j$-th coordinate and adds $a_{j-1}$ to the $(j-1)$-st coordinate exactly $k_j$ times. This implies that $\s s = \bm \psi(\s t) = \s t + A \s k$.

\medskip

Observe that the matrix $A$ has rank $m-1$, and its kernel is spanned by the vector
$$\left(\frac 1{b_{1}b_{2} \cdots b_{{m-2}} b_{{m-1}}},\frac 1{a_{1}b_{2} \cdots b_{{m-2}}b_{{m-1}}},\ldots, \frac 1{a_{1}a_{2} \cdots a_{{m-2}} a_{{m-1}}}\right) = \left( \frac 1{c_1},\frac 1 {c_2}, \ldots, \frac 1 {c_m} \right).$$
The integer vectors in $\ker A$ are integer multiples of the vector
$$\lcm(c_1,c_2,\ldots,c_m) \cdot \left( \frac 1 {c_1}, \, \frac 1{c_2}, \, \ldots \,, \, \frac 1 {c_m} \right).$$
Suppose an integer vector $\s n = (n_1,\ldots,n_m)$ with at least one negative and at least one non-negative coordinate is a solution of
\begin{equation} \label{proofs3}
 \s s = \s t + A \hskip .04cm \s n.
\end{equation}
Suppose $n_j < 0$ and $n_{j+1} \geq 0$. Then we have $-b_{j} n_j + a_{j} n_{j+1} \geq b_{j}$ and
$$b_{j} > s_j = t_j - b_{j} n_j + a_{j} n_{j+1} \geq b_{j},$$
a contradiction. This implies that the coordinates of every solution to~\eqref{proofs3} are either all negative or all non-negative. In particular, there is exactly one solution $\s n$ with
$$0 \leq n_j \leq \frac{\lcm(c_1,\ldots,c_m)}{c_j} - 1$$
for each $j$.

\medskip

First, let us show that this solution is equal to $\s k$. Denote by $n'_j$ the number of times we subtract $b_j$ from the $j$-th coordinate and add $a_{j-1}$ to the $(j-1)$-th coordinate in the first $|\s n|= n_1+\ldots+n_m$ steps of the algorithm with input $\s t$. Let $\s s' = \s t + A \hskip .004cm \s n'$. Then $\s {s'} = \s s + A (\s n' - \s n)$, and the same argument as above shows that we cannot have $n'_j - n_j < 0$, $n'_{j+1} - n_{j+1} \geq 0$. Clearly the sum of the coordinates of $\s n' - \s n$ is $0$, so we must have $\s n' = \s n$ and $\s {s'} = \s s$. This implies that $\s n = \s k$. Therefore, the maximum number of steps $\s L_\varphi(\s t)$ is at most
\begin{equation} \label{proofs2}
 \lcm(c_1,\ldots,c_m) \cdot \left( \frac 1 {c_1} + \ldots + \frac 1 {c_m} \right) - m.
\end{equation}
Consider the vector $(a_1-1,\ldots,a_m-1)$ and observe that
$$i_1 (a_1-1) +\ldots + i_m (a_m-1) = i_1 (b_1-1) +\ldots + i_m (b_m-1).$$
Since $\bm \psi(\s t)$ is an integer vector, we have $\s s = \bm \psi(\s t) = (b_1-1,\ldots,b_m-1)$. Furthermore, the vector $(-1,-1,\ldots,-1)$ is a solution of the system
$$b_j - 1 = a_j - 1 - b_j n_j + a_j n_{j+1}, \quad 1 \leq j \leq m.$$
Therefore,
$$k_j = \frac{\lcm(c_1,\ldots,c_m)}{c_j} - 1$$
and the maximum number of steps $\s L_\varphi(\s t)$ over all $\s t$, is given by the equation~\eqref{proofs2}. This finishes the proof of Theorem~\ref{main1}. \qed

\bigskip \subsection{Proof of Theorem~\ref{defs2}} \label{proofohara}

In this subsection, we present an alternative proof of O'Hara's theorem by deducing it from Theorem~\ref{main1}. Take a partition $\lambda \in \p A_n$. Without loss of generality we may assume that its support is contained in one of the connected components of $G_\varphi$. If this connected component is a cycle, the fact that O'Hara's algorithm stops and that the final result is independent of the choices made in the algorithm follows from Theorem~\ref{main1} part~\eqref{main1a}. If the connected component is a path, we can assume that $\lambda = i_1^{t_1}\cdots i_m^{t_m}$, where $i_1 \to i_2 \to \ldots \to i_m \to i_{m+1} \to \ldots$ is a part of $G_\varphi$. Assume that $i_j b_j = i_{j+1}a_{j+1}$. Remove copies of $i_1$ until their number is smaller than $b_1$. Remove copies of $i_2$ until their number is smaller than $b_2$, etc. Since $i_j$ are all different, we have $i_k > n$ for some $k$. The algorithm cannot increase the number of copies of $i_k$, so it stops before that. Furthermore, it is clear that if $k_j$ is the number of copies of $i_j$ removed, then $k_1 = \lfloor t_1/b_1 \rfloor$, $k_2 = \lfloor (t_2 + k_1 a_2)/b_2 \rfloor$, etc. Therefore the final result is independent of the order in which we remove and add parts. \qed

\bigskip \subsection{Proof of Lemma~\ref{main2}} \label{proofscxity}

Assume that the path is $\ldots \to i_{-1} \to i_0 \to i_1 \to i_2 \to \ldots$ (it can be finite, or infinite in one or both direction). A partition $\lambda$ has only a finite number of parts, without loss of generality we can assume that $\lambda = i_1^{t_1} i_2^{t_2} \cdots i_m^{t_m}$ and $\psi(\lambda) = i_1^{s_1} i_2^{s_2} \cdots i_m^{s_m}$. O'Hara's algorithm is straightforward in this case: first remove $b_1 k_1$ copies of $i_1$, where $k_1$ is the smallest integer with $t_1-b_1 k_1 < b_1$, and add $a_2k_1$ copies of $i_2$; then remove $b_2 k_2$ copies of $i_2$, where $k_2$ is the smallest integer with $t_2 + a_2 k_1 - b_2 k_2 < b_2$, and add $a_3k_2$ copies of $i_3$, etc. Clearly,
$$k_1 \leq \frac{t_1}{b_{1}} = \frac{i_1t_1}{i_1b_{1}},$$
$$k_2 \leq \frac{t_2 + a_{2} k_1}{b_{2}} = \frac{i_2t_2 + i_1b_{1} k_1}{i_2b_{2}} \leq \frac{i_1 t_1 + i_2 t_2}{i_2b_{2}},$$
$$k_3 \leq \frac{t_3 + a_{3} k_2}{b_{3}} = \frac{i_3t_3 + i_2b_{2} k_2}{i_3b_{3}} \leq \frac{i_1 t_1 + i_2 t_2 + i_3 t_3}{i_3b_{3}},$$
etc. Therefore, the total number of steps is at most
$$(i_1 t_1 + \ldots + i_m t_m) \left( \frac 1 {i_1b_{1}} + \ldots + \frac 1 {i_mb_{m}} \right).$$
Since $i_1 t_1 + \ldots + i_m t_m$ is the size $n$ of $\lambda$, the number of steps is at most $Dn$ when $D < \infty$. In general, since $i_1,\ldots,i_m$ are distinct integers, there is a $j$ with $i_j \geq m$. On the other hand, when $i_j > n$, O'Hara's algorithm stops at $i_j$. This implies that $m \leq n$, and
$$ \frac 1 {i_1b_{1}} + \ldots + \frac 1 {i_mb_{m}} \leq \frac 1 {i_1} + \ldots + \frac 1 {i_m} \leq 1 + \frac 1 2 + \ldots + \frac 1 m \leq \log m + 1 \leq \log n + 1.$$
This completes the proof. \qed

\bigskip \subsection{Proof of Theorem~\ref{main3}}

(1) For a $1$-cycle $\set i$ and a partition $\lambda$ with only parts $i$, we have $\psi(\lambda)=\lambda$ and O'Hara's algorithm takes $0$ steps. For a $2$-cycle $i \to j \to i$, by Theorem~\ref{main1} part~\eqref{main1e}, the largest number of steps of O'Hara's algorithm is equal to
$$\lcm(a_{i},b_{i})\cdot \left( \frac 1{a_{i}} + \frac 1{b_{i}}\right) - 2\leq a_{i}b_{i}\cdot \left( \frac 1{a_{i}} + \frac 1{b_{i}}\right) - 2 = a_{i} + b_{i} - 2.$$
Denote by $\lambda_{\max} = i^{a_{i}-1}j^{a_{j}-1}$ the largest partition in $\p A$. Then O'Hara's algorithm takes at most
$$a_{i} + b_{i} - 2 \leq i(a_{i}-1) + i(b_{i}-1) = |\lambda_{\max}| \hskip 0.5cm \mbox{ steps.}$$
For a partition $\lambda \vdash n \geq |\lambda_{\max}|/2$, we need, by Theorem~\ref{main1} part~\eqref{main1d}, at most $|\lambda_{\max}|\leq 2 n$ steps. Similarly, for a partition $\lambda = i^t j^s \vdash n < |\lambda_{\max}|/2$, we have $L_\varphi(\lambda) \leq \max\set{t,s} \leq n$, since $i t + j s$ is less than $i(b_{i}-1)$ or $j(b_{j}-1)$.

\medskip

Consider a cycle $C: \,i_1 \to i_m \to i_{m-1} \to \ldots \to i_1$ in $G_\varphi$ of length $m>2$. Denote by $M_C$ the maximum number of steps O'Hara's algorithm can take on partitions with support in $\set{i_1,\ldots,i_m}$. Recall that $M_C$ is given by Theorem~\ref{main1} part~\eqref{main1e}. Finally, denote by $M$ the sum of $M_C$ over all cycles $C$ of length $>2$ in $G_\varphi$.

\medskip

Note that every partition $\lambda \in \p A$ is decomposed into partitions $\lambda^{(1)},\ldots,\lambda^{(r)}$ of sizes $n_1,\ldots,n_r$ with support in one of the components of $G_\varphi$. From above and by Lemma~\ref{main2}, the number of steps of O'Hara's algorithm is at most
$$n_1 (\log n_1 + 1) + \ldots + n_r (\log n_r + 1) + 2n + M \leq n (\log n + 3) + M \leq (1+\varepsilon) n \log n,$$
for every fixed $\varepsilon > 0$ and $n$ large enough.

\medskip

(2) Take a cycle $i_1 \to i_m \to i_{m-1} \to \ldots \to i_1$ of length $m$. In the notation of Theorem~\ref{main1} part~\eqref{main1e}, we have
$$\aligned
L_\varphi(\lambda) & \, \leq \, \lcm(c_1,\ldots,c_m) \cdot \left( \frac 1 {c_1} + \ldots + \frac 1 {c_m} \right) \, - \, m \\
& \, \leq \, a_{1}a_{2}\cdots a_{{m-1}}b_{1}b_{2}\cdots b_{{m-1}} \cdot \left( \frac 1 {c_1} + \ldots + \frac 1 {c_m} \right) \, - \, m \\
& \, \leq \, \sum_{j=1}^m b_{1}b_{2}\cdots b_{{j-1}}a_{j}\cdots a_{{m-2}}a_{{m-1}} \, - \, m.
\endaligned
$$
Denote by $N$ is the largest of $\set{a_1,\ldots,a_m,b_1,\ldots,b_m}$. Then the size of the largest partition $\lambda_{\max} \in \p A$ is at least $N-1$, and by the calculation above the number of steps is at most $m (N^{m-1}-1)$. Therefore $L_\varphi(\lambda_{\max}) \leq m|\lambda_{\max}|^{m-1}$. For $\lambda \vdash n \geq |\lambda_{\max}|/m$, we have $L_\varphi(\lambda) \leq m|\lambda_{\max}|^{m-1} \leq m^m n^{m-1}$. On the other hand, when $n < |\lambda_{\max}|/m$, we have $n \leq i_j (b_{j}-1)$ for some $j$, and O'Hara's algorithm does not remove any copies of $i_j$. By Lemma~\ref{main2}, this implies that $L_\varphi(\lambda) \leq n(1+1/2+\ldots+1/m)$.

\medskip

In the notation of part (1), write $M'$ for the sum of $M_C$ over all cycles of length $> m$. As before, every partition $\lambda \in \p A_n$ is decomposed into partitions $\lambda^{(1)},\ldots,\lambda^{(r)}$ of sizes $n_1,\ldots,n_r$ with support in one of the components of $G_\varphi$. From above and by Lemma~\ref{main2}, the number of steps of O'Hara's algorithm is at most
$$m^m n_1^{m-1} + \ldots + m^m n_r^{m-1} + M' \leq m^m n^{m-1} + M' = O(n^{m-1}). \eqno \qed$$

\bigskip \subsection{Proof of Theorem~\ref{main4}} \label{proofmain4} We begin with the following key example.

\begin{exm} \label{main5}
 Let $p_1,p_2,\ldots,p_m$ be distinct primes. Take $\varphi(p_j)=p_{j+1}$ (where $p_{m+1}=p_1$), $a_{j} = p_{j+1}$, $b_{j} = p_{j-1}$ (where $p_0 = p_m$). Then
 $$a_{1}a_{2}\cdots a_{{j-1}}b_{{j}}\cdots b_{{m-2}}b_{{m-1}} = p_2 p_3 \cdots p_{j-1}p_{j}p_{j-1}p_{j}\cdots p_{m-2}$$
 for $j=1,\ldots,m$, and the lowest common multiple of these numbers is
 $$p_1p_2^2p_3^2\cdots p_{m-2}^2p_{m-1}p_m.$$
 By Theorem~\ref{main1} part~\eqref{main1e}, the maximum number of steps of O'Hara's algorithm is equal to
 $$p_1 p_2 \cdots p_m \cdot \left( \frac 1{p_1 p_2} + \frac 1{p_2 p_3} + \ldots + \frac 1{p_{m-1} p_m} + \frac 1 {p_m p_1}\right) - m.$$
 Note that Example~\ref{defs3} is a special case of this for $m=3$, $p_1=3$, $p_2=4$ and $p_3=5$.
 
 \smallskip
 
 When $p_1,\ldots,p_m$ are approximately equal, i.e.\hspace{-0.07cm} $|p_j/N - 1| < \varepsilon$ for some~$N$, we obtain partitions of relatively small size (around $m N^2$), for which O'Hara's algorithm takes a large number of steps (around $m N^{m-2}$).
\qed \end{exm}

Now we are ready to prove Theorem~\ref{main4}.

\medskip

(1) Take
$$\begin{array}{c|ccccccccccccccccccccc} i_j & 1 & \to & 2 & \to & 3 & \to & 6 & \to & 5 & \to & 10 & \to & 7 & \to & 14  & \to & 9 & \to & \ldots \\ \hline a_{j} & \infty & & 1 & & 2 & & 1 & & 6 & & 1 & & 10 & & 1 & & 14 & & \ldots \\ b_{j} & 2 & & 3 & & 2 & & 5 & & 2 & & 7 & & 2 & & 9 & & 2 & & \ldots \end{array}$$
i.e.\hspace{-0.07cm} $i_{2j-1} = 2j-1$, $i_{2j} = 4j-2$, $\varphi(i_j) = i_{j-1}$, $a_1 = \infty$, $a_{2j-1} = 4j-6$ for $j \geq 2$, $a_{4j-2} = 1$, $b_{2j-1} = 2$, $b_{4j-2} = 2j + 1$. O'Hara's algorithm on the partition $\lambda^{(k)} = 1^{2^k \cdot (2k-1)!!}$ runs as follows:
$$\begin{array}{cccl}
 \mathbf{1^{2^k \cdot 3 \cdot 5 \cdots (2k-1)}} & \to & \scriptstyle 2^{2^{k-1}\cdot 3 \cdot 5 \cdots (2k-1)} & \mbox{   in } 2^{k-1}\cdot 3 \cdot 5 \cdot 7 \cdot 9 \cdots (2k-1) \mbox{ steps}\\
 \scriptstyle 2^{2^{k-1}\cdot 3 \cdot 5 \cdots (2k-1)} & \to & \scriptstyle 3^{2^k \cdot 1 \cdot 5 \cdot 7 \cdot \cdots (2k-1)} & \\ 
 \scriptstyle 3^{2^k \cdot 1 \cdot 5 \cdot 7 \cdot \cdots (2k-1)} & \to & \scriptstyle 6^{2^{k-1}\cdot 1 \cdot 5 \cdot 7 \cdots (2k-1)} & \mbox{   in } 2^{k-1}\cdot 1 \cdot 5 \cdot 7 \cdot 9 \cdots (2k-1) \mbox{ steps}\\
 \scriptstyle 6^{2^{k-1}\cdot 1 \cdot 5 \cdot 7 \cdots (2k-1)} & \to & \scriptstyle 5^{2^k \cdot 1 \cdot 3 \cdot 7 \cdots (2k-1)} & \\ 
 \scriptstyle 5^{2^k \cdot 1 \cdot 3 \cdot 7 \cdots (2k-1)} & \to & \scriptstyle 10^{2^{k-1}\cdot 1 \cdot 3 \cdot 7 \cdots (2k-1)} & \mbox{   in } 2^{k-1}\cdot 1 \cdot 3 \cdot 7 \cdot 9 \cdots (2k-1) \mbox{ steps, etc.}
\end{array}$$
Note that $n := |\lambda^{(k)}| = 2^k (2k-1)!! = k^{\Theta(k)}$. We conclude that
$$L_\varphi(\lambda^{(k)}) \, \geq \, 2^k \cdot (2k-1)!! \left( \frac 1 {2} +\frac 1 {6} +\frac 1 {10} + \ldots + \frac 1{2(2k-1)}\right) \, = \, \Omega(n \log \log n).$$

(2) Take $D_k$ to be a superpolynomial (for example exponential) integer function of $k$. For each $k$ large enough, we can choose $m$ distinct primes $p_1^k,\ldots,p_m^k$ between $D_k$ and $2D_k$. Choose $i_j^k = km + j$ for $j=1,\ldots,m$, $a_j^k = i_{j+1}^k p_j^k$, $b_j^k = i_{j-1}^k p_{j-1}^k$ (where indices are written cyclically). Therefore, $i_1^k \to i_m^k \to \ldots \to i_1^k$ is an $m$-cycle in $G_\varphi$. The largest partition $\lambda^{(k)} \in \p A$ with support in $\set{i_1^k,\ldots,i_m^k}$ has size at most $2(k+1)m^2 D_k$. Then
$$\aligned \lcm(c_1,\ldots,c_m) & \, = \, \lcm(p_mp_1p_2\cdots p_{m-2},p_1^2p_2\cdots p_{m-2},\ldots,p_1p_2\cdots p_{m-1}) \\ & \, = \, p_1^2p_2^2\cdots p_{m-2}^2 p_{m-1}p_m \endaligned$$
where we omitted the upper index $k$ of $p_j^k$ for simplicity. Therefore,
$$L_\varphi(\lambda^{(k)}) \geq p_1p_2\cdots p_{m-1} + p_2p_3\cdots p_m + \ldots + p_mp_1 \cdots p_{m-2} - m \geq m (D_k^{m-1} - 1).$$
On the other hand, for every $C$ and $k$ large enough, we have
$$m (D_k^{m-1}-1) \geq C (2(k+1)m^2D_k)^{m-1-\varepsilon},$$
since otherwise $D_k$ grows polynomially, contradicting our assumptions. This implies that the corresponding $\overline a$, $\overline b$ and $\varphi$ satisfy the conditions of (2).

\medskip

(3) Let $p_1,p_2,\ldots$ denote the sequence of all primes. Set $\varphi(p_j)=p_{j+1}$ for $j \neq k^k$, $\varphi(p_1)=p_1$, $\varphi(p_{k^k})=p_{(k-1)^{k-1}+1}$ for $k \geq 2$, $a_{j} = p_{j+1}$ for $j \neq k^k$, $a_{1} = p_1$, $a_{k^k} = p_{(k-1)^{k-1}+1}$ for $k \geq 2$, $b_{j} = p_{j-1}$ for $j \neq 1$, $j \neq k^k+1$, $b_{1} = p_1$, $b_{{k^k+1}}=p_{(k+1)^{k+1}}$ for $k \geq 1$. The following table summarizes these values.
$$\begin{array}{c|c|ccc|cccccccccccccc|c}
 i_j & 2 & 3 & 5 & 7 & 11 & 13 & 17 & 19 & 23 & 29 & 31 & 37 & 41 & 43 & 47 & \ldots & 101 & 103 & \ldots \\
 \hline a_{j} & 2 & 5 & 7 & 3 & 13 & 17 & 19 & 23 & 29 & 31 & 37 & 41 & 43 & 47 & 53 & \ldots & 103 & 11 & \ldots \\
 b_{j} & 2 & 7 & 3 & 5 & 103 & 11 & 13 & 17 & 19 & 23 & 29 & 31 & 37 & 41 & 43 & \ldots & 97 & 101 & \ldots
\end{array}$$
Clearly, $G_\varphi$ consists of cycles of length $1$ and $k^k-(k-1)^{k-1}$ for $k \geq 2$.
For $k \geq 2$, write $m_1=(k-1)^{k-1}$, $m_2=k^k$, and define $\lambda^{(k)} \in \p A$ to be the partition
$$(p_{m_1+1})^{(p_{m_1+2})-1}(p_{m_1+2})^{(p_{m_1+3})-1}\cdots (p_{m_2-1})^{(p_{m_2})-1}(p_{m_2})^{(p_{m_1+1})-1}$$
of size
$$n_k = p_{m_1+1} p_{m_1+2} + p_{m_1+2}p_{m_1+3} + \ldots + p_{m_2-1}p_{m_2} + p_{m_2}p_{m_1+1} - p_{m_1+1} - p_{m_1+2} - \ldots - p_{m_2}.$$
By the calculation in Example~\ref{main5}, O'Hara's algorithm takes exactly
$$p_{m_1+1}p_{m_1+2}\cdots p_{m_2-2} + p_{m_1+2}p_{m_1+3}\cdots p_{m_2-1} + \ldots + p_{m_2}p_{m_1+1}\cdots p_{m_2-3} - (m_2 - m_1)$$
steps to compute $\psi(\lambda^{(k)})$. By the distribution law of prime numbers, we have $p_n = n \log n (1 + o(1))$. Therefore,
$$n_k \leq (m_2-m_1) p_{m_2}^2 = (k^k-(k-1)^{k-1}) p_{k^k}^2 \sim k^k (k^k \log(k^k))^2 = k^{3k+2} (\log k)^2$$
and
$$\log L_\varphi(\lambda^{(k)}) \geq \log\left( (m_2-m_1)(p_{m_1+1}^{m_2-m_1-2}-1) \right) \sim k^k \log\left((k-1)^{k-1}\right) \sim k^{k+1} \log k.$$
Thus, $L_\varphi(\lambda^{(k)}) > \exp(\sqrt[3]{n_k})$ for $k$ large enough, as desired. \qed

\bigskip \subsection{Proof of Proposition~\ref{main6}} \label{proofmain6}

By construction of the algorithm, $\s k \in \Z^m$ and it satisfies the inequalities in Proposition~\ref{main6}. Therefore, it suffices to prove that $|\s k| < |\s n|$ for every non-negative integer vector $\s n \neq \s k$ satisfying $\s 0 \leq \s t + A \s n \leq \s b - \s 1$.

\medskip

Assume that $|\s n| \leq |\s k|$. Denote by $n_j'$ the number of times we remove $b_j$ copies of $i_j$ in the first $|\s n|$ steps of O'Hara's algorithm. Define $\s {s'} = \s t + A \hskip 0.04cm \s n'$. From above, $\s {s'} = \s s + A (\s n' - \s n)$. If there exists $j$ such that $n'_j - n_j > 0$ and $n'_{j+1} - n_{j+1} \leq 0$, then
$$0 \leq s_j' = s_j - b_{j} (n'_j - n_j) + a_{j} (n'_{j+1} - n_{j+1}) < b_{j} + (-b_{j})=0,$$
which is impossible. Therefore, all the coordinates of $\s n' - \s n$ have the same sign. On the other hand, the sum of the coordinates of $\s n' - \s n$ is $0$, which implies $\s n' = \s n = \s k$. Therefore, for all $\s n \neq \s k$ as above, we have $|\s n| > |\s k|$, as desired. \qed

\bigskip

\section{Application: the speedy O'Hara's algorithm} \label{speedy}

In~\cite[\S 8.2]{P3}, a simple speed-up of O'Hara's algorithm was given: at each step, if the number of parts $i$ in $\lambda \vdash n$ is $r b_i + s$ for $0 \leq s < b_i$, remove $rb_i$ copies of part $i$ and add $ra_j$ copies of part $j$, where $\varphi(j)=i$. This replaces $r$ steps of O'Hara's algorithm with one step. This algorithm produces the same output as the original algorithm, but the number of steps may depend on the choices we make in the execution. It is called the \emph{speedy O'Hara's algorithm}. Observe now that it cannot be much faster than the usual O'Hara's algorithm: since $r \leq i(r b_i + s) \leq n$, the number of steps can be reduced by at most the order of $n$, where $\lambda \vdash n$. Thus the reasoning from the proof of Theorem~\ref{main4} part~\eqref{main4c} still gives superpolynomial lower bounds for the number of steps. The following examples show that the speed-up is not substantial (for example, of logarithmic complexity) even in the case described in Theorem~\ref{main3} part~\eqref{main3a}, and that it gives the same bound as Theorem~\ref{main3} part~\eqref{main3b}.

\begin{exm}
 Take
 $$\begin{array}{c|ccccccccccccccccccc} i_j & \ldots &  \to & 18 & \to & 7 & \to & 14 & \to & 5 & \to & 10 & \to & 3 & \to & 6 & \to & 1\\ \hline a_{j} & \ldots & & 1 & & 18 & & 1 & & 14 & & 1 & & 10 & & 1 & & 6 \\ b_{j} & \ldots & & 7 & & 2 & & 5 & & 2 & & 3 & & 2 & & 1 & & \infty \end{array}$$
 i.e.\hspace{-0.07cm} $i_{2j-1} = 2j-1$, $i_{2j} = 4j+2$, $\varphi(i_j) = i_{j+1}$, $a_{2j-1} = 4j+2$, $a_{4j+2} = 1$, $b_1=\infty$, $b_{2j-1} = 2$ for $j \geq 2$, $b_{4j+2} = 2j - 1$. This is similar to the example in the proof of Theorem~\ref{main4} part~\eqref{main4a}. The speedy O'Hara's algorithm on $\lambda^{(k)} = (2k-1)^{4k-2}$ for $k \geq 5$ runs as follows:
 $$\mathbf{(2k-1)^{4k-2}}  \scriptstyle \: \to \:  (4k-2)^{2k-1} \: \to \: (4k-2)^2 (2k-3)^{4k-2} \: \to \: (4k-2)^2(4k-6)^{2k-1} \: \to \: (4k-2)^2(4k-6)^4(2k-5)^{4k-6}$$
 $$\scriptstyle \: \to \: (4k-2)^2(4k-6)^4(4k-10)^{2k-3} \: \to \: (4k-2)^2(4k-6)^4(4k-10)^4 (2k-7)^{4k-10} \: \to \: \ldots \: \to \: (4k-2)^2(4k-6)^4\cdots 14^{4} 10^{7}$$
 $${\scriptstyle \: \to \: (4k-2)^2(4k-6)^4\cdots 14^{4} 10^{1}3^{20} \: \to \: (4k-2)^2(4k-6)^4 \cdots 14^{4} 10^{1}6^{10} \: \to \:} \mathbf{(4k-2)^2(4k-6)^4 \cdots 14^{4} 10^{1}1^{60}}.$$
 In particular, the speedy O'Hara's algorithm takes $2k-2=\Theta(\sqrt n)$ steps for a partition of size $n=(2k-1)(4k-2)=2(2k-1)^2$.
\qed \end{exm}

Note that the proof that $m \leq n$ in the proof of Lemma~\ref{main2} together with the proof of Theorem~\ref{main3} part~\eqref{main3a} shows that the speedy O'Hara algorithm takes $O(n)$ steps when $G_\varphi$ has a finite number of cycles of length greater than $2$. This is smaller than $\Omega(n \log \log n)$ obtained earlier. In other words, the speedy O'Hara's algorithm can be asymptotically faster.

\begin{exm}
 Take $\overline a$, $\overline b$ and $\varphi$ constructed in the proof of part~\eqref{main4b} of Theorem~\ref{main4}, see Subsection~\ref{proofmain4}. For each $k$, we have $i_m^{k}/i_1^k < (k+1)/k\leq 2$, so at each step the number of parts $i_j$ is at most
 $$\frac 1 {i_j^k} \left({i_1^k (i_2^k p_1^k - 1) + \ldots + i_m^k (i_1^k p_m^k - 1)}\right) \leq 8m i_1^k p_1^k.$$
 Therefore, in the speedy algorithm we replace $b_j^k=i_{j-1}^k p_{j-1}^k$ copies of $i_j$ by $a_j^k$ copies of $i_{j-1}$ at most $8m$ times. This means a speed-up by only a constant factor, and the number of steps of the speedy O'Hara's algorithm is $\Omega(n^{m-1-\varepsilon})$ for every $\varepsilon > 0$.
\qed \end{exm}

\bigskip

\section{Final remarks} \label{final}

 \bigskip

\subsection{} \label{final1}
The polynomial time
algorithm in the proof of Theorem~\ref{linprog1} is given implicitly, by using the general results in integer linear programming. It is saying that the function $\psi: \mathcal A_n \to \mathcal B_n$ can be computed much faster, by circumventing the elegant construction of O'Hara's algorithm. It would be interesting to give an explicit construction of such an algorithm.

\medskip

In a different direction, it might prove useful to restate other involution principle bijections in the language of linear programming, such as the Rogers-Ramanujan bijection in~\cite{GM2} or in~\cite{BoP}. If this works, this might lead to a new type of a bijection between these two classes of partitions. Alternatively, this might resolve the conjecture by the second author on the mildly exponential complexity of Garsia-Milne's Rogers-Ramanujan bijection, see \cite[Conjecture 8.5]{P3}.


\subsection{}
Note the gap between the number $\exp \Theta(\sqrt n)$ of partitions of~$n$ and the lower
bound $\p L_\varphi(n) = \exp \Omega(\sqrt[3]{n})$ in Theorem~\ref{main4}.
It would be interesting to decide which of the two worst complexity bounds on the number of steps of O'Hara's algorithm is closer to the truth.

\medskip

Note that we applied our linear programming approach only in the bounded cycle case. We do not know if there is a way to apply the same technique to the general case. However, we believe that there are number theoretic obstacles preventing that and in fact, computing O'Hara's bijection as a function on partitions  may be hard in the formal complexity sense.


\subsection{}
It would be interesting to find~$\varphi$ such that the graph~$G_\varphi$ is a path
and $\p L_\varphi(n) = \Theta(n \log n)$. From the proof of Lemma~\ref{main2}
it follows that the number of steps of (the usual) O'Hara's algorithm is at
most $n (\log k + 1)$,  where~$k$ is the number of steps of the speedy O'Hara's algorithm.
In Subsection~\ref{speedy} we constructed an example with $k \sim \sqrt {n/2}$,
so it certainly seems possible that such examples exist.


\subsection{}
Most recently, variations on the O'Hara's bijection and applications of rewrite systems were found in~\cite{SSM} and~\cite{K1,K2}. It would be interesting to see connections between our analysis and this work.


\subsection{}
In the finite dimensional case, the structure of the map~$\mathbf{\psi}$ establishing
approximate $\Pi$-congruence remains largely unexplored. For example,
it would be nice to obtain some convergence result in the
irrational case using the rational approximations which follow from
part~$(3)$ of Theorem~\ref{main1} and upper bounds in part~$(2)$ of
Theorem~\ref{main3}.

\medskip

Recall also that the $2$-dimensional case can be viewed as the Euclid algorithm
which in turn corresponds to the usual continued fractions (see Example~\ref{ex-euclid}).
Thus the geometry of~$\mathbf{\psi}$  can be viewed as a delicate multidimensional
extension of continued fractions. Given the wide variety of (different) multidimensional
continued fractions available in the literature, it would be interesting to see
if there is a connection to at least one of these notions.

\vskip 0.5cm

{\bf Acknowledgments.} \  We are grateful to George Andrews and Dennis Stanton
for their interest in the paper and to Kathy O'Hara for sending us a copy of her
thesis~\cite{O1}. The second named author was supported by the NSF. He would
also like to thank Vladimir Arnold, Elena Korkina and Mark Sapir for teaching
him about multidimensional continued fractions.

\vskip 1cm

\bigskip

\bigskip

\bigskip

{\sc \scriptsize Department of Mathematics, Massachusetts Institute of Technology, Cambridge, MA 02139\\
\tt{konvalinka@math.mit.edu}\\
\tt{http://www-math.mit.edu/\~{}konvalinka/}}

\bigskip

{{\sc \scriptsize School of Mathematics, University of Minnesota, Minneapolis, MN 55455}{\hskip 0.5cm \scriptsize and}\\
\sc \scriptsize Department of Mathematics, Massachusetts Institute of Technology, Cambridge, MA 02139\\
\tt{pak@math.mit.edu}\\
\tt{http://www-math.mit.edu/\~{}pak/}}


\begin{thebibliography}{123456}

\bibitem[A]{A}  G.~E. Andrews,
\emph{The theory of partitions} (Second ed.),
Cambridge U. Press, Cambridge, 1998.

\bibitem[BP]{BoP}
C. Boulet and I. Pak,
A combinatorial proof of the Rogers-Ramanujan identities,
\emph{J. Combin. Theory Ser. A}~\textbf{113} (2006), 1019--1030.

\bibitem[GM1]{GM1}
A.~M. Garsia and S.~C. Milne,
Method for constructing bijections for classical partition
identities, \emph{Proc. Nat. Acad. Sci. U.S.A.}~\textbf{78} (1981),
no.~4,  2026--2028.

\bibitem[GM2]{GM2}
A.~M. Garsia and S.~C. Milne,
A {R}ogers-{R}amanujan bijection
\emph{J. Combin. Theory Ser. A}~\textbf{31} (1981),
289--339.

\bibitem[G]{G}
B. Gordon,
Sieve-equivalence and explicit bijections,
\emph{J. Combin. Theory Ser. A}~\textbf{34} (1983),
90--93.

\bibitem[K1]{K1}
M. Kanovich,
Finding direct partition bijections by two-directional rewriting techniques,
\emph{Discrete Math.}~\textbf{285}  (2004),  151--166.


\bibitem[K2]{K2}
M. Kanovich,
The two-way rewriting in action: removing the mystery of Euler-Glaisher's map,
\emph{Discrete Math.}~\textbf{307}  (2007),  1909--1935.

\bibitem[O1]{O1}
K.~M. O'Hara,
\emph{Structure and Complexity of the Involution Principle for Partitions},
Ph.D. thesis, UC Berkeley, California, 1984, 135 pp.

\bibitem[O2]{O2}
K.~M. O'Hara,
Bijections for partition identities,
\emph{J. Combin. Theory Ser. A}~\textbf{49} (1988), 13--25.

\bibitem[P1]{P1}
I. Pak,
Partition identities and geometric bijections,
\emph{Proc. A.M.S.}~\textbf{132} (2004), 3457--3462.

\bibitem[P2]{P2}
I. Pak,  The nature of partition bijections I. Involutions,
\emph{Adv. Applied Math.}~\textbf{33} (2004), 263--289.

\bibitem[P3]{P3}
I. Pak, Partition bijections, a survey, \emph{Ramanujan J.}~12 (2006),
5--75.

\bibitem[P4]{P4} I. Pak,
The nature of partition bijections~II. Asymptotic stability,
preprint, 32 pp., available at {\tt http://www-math.mit.edu/\~{}pak/}


\bibitem[PV]{PV}
I. Pak and E. Vallejo,
Combinatorics and geometry of Littlewood-Richardson cones,
\emph{Europ. J. Combin.}~\textbf{26} (2005), 995--1008.

\bibitem[R]{R}
J.~B. Remmel,
Bijective proofs of some classical partition identities.
\emph{J. Combin. Theory Ser. A}, \textbf{33} (1982), 273--286.

\bibitem[S]{S}
A. Schrijver, \emph{Theory of linear and integer programming}, John Wiley,
Chichester, 1986.

\bibitem[SSM]{SSM}
J.~A. Sellers,  A.~V. Sills and G.~L. Mullen,
Bijections and congruences for generalizations of partition identities of Euler and Guy,
\emph{Electronic J. Combin.}~ \textbf{11} (2004), no.~1, RP~43, 19 pp.

\end{thebibliography}
\end{document}